\documentclass[11pt]{amsart}
\usepackage{amsmath,amscd}
\title[Mordell-Weil Representations]{Realizing Rational Representations in Mordell-Weil Groups}  
\author{Bo-Hae Im and Michael Larsen}
\date{February 1, 2010}
\address{Department of Mathematics, Chung-Ang University, 221, Heukseok-dong, Dongjak-gu, Seoul, 155-756, South Korea}\email{bohaeim@gmail.com}
\address{Department of Mathematics, Indiana University, Bloomington,
Indiana 47405, USA} \email{larsen@math.indiana.edu}
\thanks{Michael Larsen was partially supported by NSF grants DMS-0100537 and  DMS-0800705.}

\newtheorem{thm}{Theorem}[section]
\newtheorem{cor}[thm]{Corollary}
\newtheorem{prop}[thm]{Proposition}

\theoremstyle{remark}

\theoremstyle{definition}

\newcommand\g{{\bf g}}
\newcommand\p{{\bf p}}
\renewcommand\P{{\mathbb P}}

\newcommand{\F}{{\mathbb F}}

\newcommand{\Z}{{\mathbb Z}}

\newcommand{\C}{{\mathbb C}}
\newcommand{\Q}{{\mathbb Q}}

\newcommand{\Lie}{\operatorname{Lie}}

\newcommand{\SL}{{\rm SL}}

\newcommand{\SO}{{\rm SO}}
\newcommand{\Gal}{{\rm Gal}}

\newcommand{\sing}{\operatorname{sing}}
\newcommand{\tor}{\operatorname{tor}}

\newcommand{\Hom}{\operatorname{Hom}}

\newcommand{\Aut}{\operatorname{Aut}}
\newcommand{\Ind}{\operatorname{Ind}}
\newcommand{\Spec}{\operatorname{Spec}\,}

\begin{document}

\begin{abstract}
Let $G$ be a finite group and $V$ a finite-dimensional rational $G$-representation.
We ask whether there exists a finite Galois extension $L/K$ of number fields
with Galois group $G$, an elliptic curve $E/K$, and a $G$-submodule of $E(L)\otimes \Q$
isomorphic to $V$.  
\end{abstract}

\maketitle

\section{Introduction}\label{intro}

Let $L/K$ be a Galois extension of number fields with group $G$.  Let $A/K$ be an abelian variety.  We regard $A(L)\otimes\Q$ as a $\Q[G]$-module.  We say that a faithful, irreducible, rational representation $V$ of a finite group $G$ is \emph{Mordell-Weil}, or more specifically, \emph{Mordell-Weil in dimension $\dim A$} or \emph{Mordell-Weil over $K$ in dimension $\dim A$}, 
if $V$ arises as a $G$-submodule of $A(L)\otimes\Q$ for some triple $(K,L,A)$.  It turns out that every rational irreducible representation of a finite group 
is Mordell-Weil:

\begin{thm}
\label{rational}
If $G$ is a finite group, $V$ is a finite-dimensional $\Q$-vector space and $\rho:G\hookrightarrow \Aut(V)$ is a faithful irreducible rational representation, 
then there exists a finite Galois extension of number fields, $L/K$, 
an isomorphism $\Gal(L/K)\to G$, and a $\Gal(L/K)$-stable subspace $W$ of
$A(L)\otimes\Q$ such that by transport of structure, $W$ is isomorphic to $V$ as a 
$\Q[G]$-module.
\end{thm}

We are particularly interested in the case $\dim A = 1$, and we present two methods, one geometric, and one arithmetic, for showing that certain interesting pairs $(G,V)$
can be realized inside the Mordell-Weil groups of elliptic curves.  These two methods are
loosely analogous to two established approaches to the inverse Galois problem, namely the rigidity method, and the use of automorphic forms.

\section{Hilbert Irreducibility}

The geometric method of realizing a pair $(G,V)$ was inspired by \cite{Im}.  The idea is to look for 
a diagram
\begin{equation*}
\begin{CD}
X @>>> A \\
@VVV  @VVV \\
\P^1 @>>> A/G \\
\end{CD}
\end{equation*}
over a number field $K$.  Here, the first row consists of a 
$G$-equivariant map from a curve to an abelian variety and
the vertical arrows are quotient maps.
Hilbert irreducibility gives a $K$-point in $\P^1 = X/G$ whose preimage in $X$ consists of a single point $\Spec L$, where $G=\Gal(L/K)$.  This point generates a 
$\Q[G]$-submodule of $A(L)\otimes\Q$; the only difficulty is to insure that this module contains $V$
as a submodule.   This can be achieved by the following proposition:

\begin{prop}
\label{realize}
Let $A/K$ and $B/K$ be abelian varieties and $G$ a finite group of $K$-automorphisms of $A$.
Suppose $A$ contains an irreducible curve $X$ which is stabilized by $G$ and is not contained
in the translate of any proper abelian subvariety of $A$.   Let $V$ denote any $\Q$-irreducible subrepresentation of $\Hom_K(A,B)\otimes\Q$.  If $X/G$ is a rational curve over $K$, then there exists a Galois extension $L/K$ with group $G$ and a $\Q[G]$-submodule of $B(L)\otimes\Q$ isomorphic to $V$.  In particular, $(G,V)$ is Mordell-Weil.
\end{prop}

\begin{proof}
Let $v$ denote any non-zero vector in $V$ and $\phi\in\Hom_K(A,B)$ a non-zero scalar multiple of $v$.  
Let $P = X/G$ and $\pi$ the quotient morphism $X\to P$.
For $g\neq 1$, $X$ is not contained in the kernel of $1-g$ acting on $A$.  Therefore,
the morphism $X\to P$ is a regular branched covering with Galois group $G$.
By an observation of Silverman \cite{Si}, originally made in the setting of elliptic curves, but equally valid for abelian varieties, the set of torsion points on $B$ which are defined over number fields of degree $\le |G|$ over $K$ is finite.  Therefore, the set of $p\in P(K)$ such that 
$\pi^{-1}(p)\cap \phi^{-1}B(\bar K)_{\tor}\neq\emptyset$ is finite.  By the Hilbert irreducibility theorem, there exists $x\in X(\bar K)$ such that $\pi(x)\in P(K)$, 
the $\Gal(\bar K/K)$-orbit of $x$ coincides with the $G$-orbit of $x$, and $\phi(x)$ is a point of infinite order on $B$.

Evaluation at $x$ gives a $\Q[G]$-linear map $\Hom_K(A,B)\otimes\Q\to B(L)\otimes\Q$.  
As composition with the inclusion $V\hookrightarrow\Hom_K(A,B)\otimes\Q$ is non-zero, $B(L)\otimes\Q$ contains at least one copy of $V$.
\end{proof}

Given a finite group $G$, an $n$-tuple $\g = (g_1,\ldots,g_n)\in G^n$ satisfying 
$$g_1 g_2\cdots g_n = 1,$$
and an $n$-tuple 
$\p = (p_1,\ldots,p_n)\in \P^1(\overline{\Q})^n$ whose coordinates are pairwise distinct, we define $U_{\g,\p}/\C$ as the open curve such that $U_{\g,\p}(\C)$ is the regular covering space of $\P^1(\C)\setminus\{p_1,\ldots,p_n\}$ with deck transformations in $\langle g_1,\ldots,g_n\rangle$ and local monodromy $g_i$ at $p_i$. 
We assume $G = \langle g_1,\ldots,g_n\rangle$, so $G$ acts freely on $U_{\g,\p}$.  
Let $X_{\g,\p}$ denote the non-singular compactification of $U_{\g,\p}$.  The action of $G$ on $U_{\g,\p}$ extends to $X_{\g,\p}$, and $\P^1 = X_{\g,\p}/G$.  Moreover, $X_{\g,\p}$ can be defined over some number field.

By a theorem of Weil \cite[VI~Prop.~7]{Se},
 the representation of $G$ on $H^1(X_{\g,\p}(\C),\Q)$ satisfies
\begin{equation}
\label{RH}
H^1(X_{\g,\p}(\C),\Q) \oplus I_G \oplus I_G \cong \bigoplus_{i=1}^n \Ind_{\langle g_i\rangle}^G I_{\langle g_i\rangle}.
\end{equation}
An immediate consequence of this is that 
for any complex representation $V_{\C}$ of $G$,
$$2 \dim V_{\C} - 2\dim V_{\C}^G  \le \sum_{i=1}^n \dim V_{\C} - \dim V_{\C}^{g_i}.$$
This is the finite case of a well-known result of L.~Scott \cite{Sc}.
We are interested in the case $V_{\C} = V\otimes_{\Q}\C$, where $V$ is a $\Q$-vector space 
on which $G$ acts.  Here we have the following:

\begin{prop}
If $V$ is a rational representation, then 
\begin{equation}
\label{genus}
-2 \dim V + 2\dim V^G  + \sum_{i=1}^n \dim V - \dim V^{g_i} = 2g
\end{equation}
for some non-negative integer $g$.
\end{prop}

\begin{proof}
The space
$$W = \Hom_G(H^1_{\sing}(X_{\g,\p}(\C),\Q),V)$$
inherits a rational Hodge structure of weight 1 from $H^1_{\sing}(X_{\g,\p}(\C),\Q)$ and is therefore of even dimension; on the other hand, its dimension is given by the left hand side of (\ref{genus}).

\end{proof}

We call $g$ the \emph{genus} of the triple $(G,V,\g)$.
We can now prove a more precise version of Theorem~\ref{rational}.

\begin{thm}
Let $(G,V,\g)$ is a triple consisting of a finite group, 
an irreducible rational representation space, 
and a generating $n$-tuple with product $1$ such that $(G,V,\g)$ has genus $g>0$.  
Then $(G,V)$ is Mordell-Weil in dimension $g$, i.e.,
there exists a number field $K$, an abelian variety $A/K$ of dimension $g$, a Galois extension $L/K$, an isomorphism $\Gal(L/K)\to G$, and a $\Gal(L/K)$-stable subspace W of
$A(L)\otimes\Q$ such that by transport of structure, $W$ is isomorphic to $V$ as $\Q[G]$-module.
\end{thm}

\begin{proof}
Let $v$ denote a non-zero vector in the representation space $V$, and let $\Lambda = \Z[G]v$
denote the corresponding lattice.  
If $A/K$ is an abelian variety with a $G$-action, we define 
$A_\Lambda :=A\otimes_{\Z[G]}\Lambda$ to be the 
functor represented by $S\mapsto A(S)\otimes_{\Z[G]}\Lambda$ for every scheme $S$ over $K$.
Concretely, $A_\Lambda$ is the quotient of $A$ by $\sum \alpha A$, where the sum is taken over 
$\alpha\in\ker (\Z[G]\to\Lambda)$.  As the quotient of an abelian variety by a closed subgroup, 
$A_\Lambda$ is again an abelian variety.  
Its Lie algebra is $\Lie(A)\otimes_{\Q[G]}V_\Q$.  The vector space $\Hom_K(A,A_\Lambda)\otimes\Q$ admits a $G$-action (given by the action of $G$ on $A$) and contains a non-zero vector $e$ (the natural quotient map) provided that $\Lie(A)\otimes_{\Q[G]}V_\Q$, and therefore $A_\Lambda$, is non-zero.   On the other hand,
$e$ is annihilated by $\ker (\Z[G]\to\Lambda)$.  It follows that $\Q[G]e\cong V$ as $\Q[G]$-module.

Fix $\p\in\P^1(\overline{\Q})^n$.  Let $D$ denote the divisor of $X_{\g,\p}$ which is the inverse image under the map $X_{\g,\p}$ of the divisor $[0]$ on $\P^1$.
Let $J$ denote the Jacobian variety of $X_{\g,\p}$. 
By hypothesis, 
\begin{equation*}
\begin{split}
\dim J_\Lambda = \frac{\dim \Lie(J)\otimes_{\Q[G]}V_\Q}{2} 
&= \frac{\dim H^1_{\sing}(X_{\g,\p}(\C),\Q)^*\otimes_{\Q[G]}V}{2} \\
&= \frac{\dim H^1_{\sing}(X_{\g,\p}(\C),\Q)\otimes_{\Q[G]}V}{2}= g > 0. \\
\end{split}
\end{equation*}

The morphism $X_{\g,\p}\to J$ given by $Q\mapsto |G|[Q]-D$ is $G$-equivariant, and the quotient $X_{\g,\p}/G$ is isomorphic to $\P^1$.   The theorem follows by applying Proposition~\ref{realize} to the morphism $J\to J_\Lambda$. 

\end{proof}

We remark that the Hilbert irreducibility argument actually gives a little more: it shows that we may choose infinitely many linearly disjoint extensions $L_i$ over $K$, all with Galois group $G$ and submodules $V_i\subset A(L_i)$ which are isomorphic to $V$ as $\Gal(L_i/K)=G$-modules.

\section{Some genus $1$ triples}

\begin{prop}
\label{Weyl}
Let $G$ be the Weyl group of an irreducible root system of rank $r$ 
and $V$ its reflection representation.  
Let $G^\circ = G\cap \SO(V)$.
Then there exist vectors $\g\in G^{2r+2}$ and $\g^\circ\in (G^\circ)^{r+1}$ such that 
$(G,V,\g)$ and $(G^\circ,V,\g^\circ)$ have genus 1.
\end{prop}

\begin{proof}
Let $s_1,\ldots,s_r$ denote the simple reflections.  We can take
$$\g = (s_1,s_1,s_2,s_2,\ldots,s_r,s_r,s_1,s_1).$$
As $\dim V^{s_i} = r-1$, we have $g=1$.
Every Dynkin diagram with $r-1\ge 2$ edges can be written as the union of two paths which meet at
a single vertex:
$i_1,\ldots,i_p$ and $j_1,\ldots,j_q$, with $p+q = r+1$.   Then we can take
$$\g^\circ = (s_{i_1}s_{i_2},s_{i_2}s_{i_3},\ldots,s_{i_p}s_{i_1},
s_{j_1}s_{j_2},s_{j_2}s_{j_3},\ldots,s_{j_q}s_{j_1}).$$
Every product of the form
$s_{i_m}s_{i_n}$ or the form $s_{j_m}s_{j_n}$ is obviously in the group $\langle \g^\circ\rangle$
generated by the coordinates of $\g^\circ$, and as there exists one vertex of the form $i_k=j_l$,
every product $s_{i_m}s_{j_n} = s_{i_m}s_{i_k}s_{j_l}s_{j_n}$ is again in $\langle \g^\circ\rangle$.
It follows that this group is the kernel of the determinant map on $G$.  Each product of two simple reflections fixes a subspace of $V$ of codimension 2, so again $g=1$.
\end{proof}

\begin{prop}
\label{Conway}
If $G$ is the automorphism group $\textup{2.Co}_1$ of the Leech lattice $\Lambda_{24}$ and $V = \Lambda_{24}\otimes\Q$, then there exists $\g\in G^3$ such that $(G,V,\g)$ has genus 1.
\end{prop}

\begin{proof}
By \cite{DA}, $\textup{Co}_1$ has a $(2A,7B,13A)$ generation, and this can be lifted to
a generating triple $(\widetilde{\text{2A}},\widetilde{\text{7B}},-\widetilde{\text{13A}})$
in $\textup{2.Co}_1$, where the lifts $\widetilde{\text{7B}}$ and $\widetilde{\text{13A}}$ are chosen 
to have order $7$ and $13$ respectively; the lift of 2A is then determined by the product $1$ condition.  By \cite{Atlas}, the resulting triple has genus 1.

\end{proof}

There is an extensive literature devoted to pairs of elements $(g_1,g_2)$ generating sporadic simple groups $G$.   In particular, cases in which the orders of $g_1$, $g_2$, and $g_1 g_2$ are all low have been extensively studied, in an attempt to classify simple Hurwitz groups and, more generally, to compute the symmetric genera of sporadic groups.  Any generating pair $(g_1,g_2)$ gives rise to a triple $\g = (g_1,g_2,g_2^{-1}g_1^{-1})$ as above.
The following table, giving some examples of genus 1, is mainly extracted from this literature.

\medskip

\begin{tabular}{|r|c|c|c|c|l|}\hline
Group & Character & Dimension & Generators & Reference \\ \hline
M$_{11}$&$\chi_2$&10&see (M$_{11}$)&\\
M$_{12}$&$\chi_2$&11&see (M$_{12}$)&\\
M$_{22}$&$\chi_2$&21&see (M$_{22}$)&\\
M$_{23}$&$\chi_2$&22&see (M$_{23}$)&\\
HS&$\chi_2$&22&2B, 5B, 7A&\cite{GM2}\\
McL&$\chi_2$&22&2A, 5A, 8A&\cite{CWW}\\
M$_{24}$&$\chi_2$&23&see (M$_{24}$)&\\
Co$_3$&$\chi_2$&23&2B, 3C, 11A&\cite{GM1}\\
Co$_2$&$\chi_2$&23&2B, 5A, 11A&\cite{GM3}\\
2.Co$_1$&$\chi_{102}$&24&$\widetilde{\text{2A}}$, $\widetilde{\text{7B}}$, $-\widetilde{\text{13A}}$&\cite{DA}\\
Tits&$\chi_6$&78&2A, 3A, 13A&\cite{AI}\\
J$_2$&$\chi_{12}$&160&2B, 3B, 7A&\cite{Wo}\\
\hline

\end{tabular}
\medskip

Where no reference is given, the assertions can easily be checked by machine, e.g., using
\cite{GAP}.  Using \cite{Atlas} notation (except that for the large Mathieu groups we write $A,B,C,\ldots,X$ instead of $0,1,\ldots,22,\infty$), we have generating pairs as follows:
\begin{gather*}
\tag{M$_{11}$}
(0183649X257)\ (07365481)(29) = (2X)(34)(59)(67) \\
\tag{M$_{12}$}
(058263X4179)\ (0\infty92)(13)(458X)(67) = (0\infty)(1X)(25)(37)(48)(69) \\
\tag{M$_{22}$}
\begin{split}
\mathtt{(AFMIHBLCRPD)(EGSVQJNOUKT)}\ &\mathtt{(ADLQF)(BHMVJ)(CPTUO)(ERNSG)}\\
= &\mathtt{(CD)(EP)(HI)(JL)(KT)(MQ)(NV)(OR)}\\
\end{split} \\
\tag{M$_{23}$}
\begin{split}
\mathtt{(AWEIHURTPBCSLGMOKJVNFD)}\ &\mathtt{(ADBPW)(CFNJS)(ETOUH)(GLKRQ)}\\
= &\mathtt{(CD)(EP)(HI)(JL)(KT)(MQ)(NV)(OR)}\\
\end{split} \\
\tag{M$_{24}$}
\begin{split}
\mathtt{(ATSX)(DW)(EQIG)}&\mathtt{(FULV)(HJOM)(NP)}\ \mathtt{(ASJVOTIFHPWBDNMLCUQKERG)}\\
= &\mathtt{(AX)(BW)(CL)(DP)(ER)(FJ)(GT)(HN)(IU)(KQ)(MV)(OS)}\\
\end{split} \\
\end{gather*}

\section{Modular curves and Steinberg representations}

If $n$ is a positive integer, we define as usual
$$\Gamma(n) := \ker \SL_2(\Z)\to \SL_2(\Z/n\Z)$$
and
$$\Gamma_0(n) :=  \Bigl\{\Bigl(\begin{matrix}a&b\\c&d\end{matrix}\Bigr)\in\SL_2(\Z)\Bigm|c\equiv 0\pmod n\Bigr)\Bigr\}.$$
Let 
$$\Gamma_{m,n} := \Gamma(m)\cap \Gamma_0(n).$$
Thus $\Gamma_{m,n}$ is normal in $\Gamma_0(n)$, and there is a natural inclusion homomorphism
$$\Gamma_0(n)/\Gamma_{m,n} \to \SL_2(\Z)/\Gamma(m) \cong \SL_2(\Z/m\Z).$$
If $m$ and $n$ are relatively prime, this inclusion is an isomorphism, by the Chinese remainder 
theorem.  Let $Y_0(n)$ and  $Y_{m,n}$ denote the quotient of the upper half-plane by
$\Gamma_0(n)$ and $\Gamma_{m,n}$ respectively.  If $(m,n)=1$,
$\SL_2(\Z/m\Z)$ acts faithfully on $Y_{m,n}$ with quotient $Y_0(n)$.
Letting $X_0(n)$ (resp. $X_{m,n}$) denote the non-singular compactification of 
$Y_0(n)$ (resp. $Y_{m,n}$), the $\SL_2(\Z/m\Z)$-action on $Y_{m,n}$ extends uniquely to 
$X_{m,n}$, and the quotient is $X_0(n)$.

\begin{thm}
\label{steinberg}
If there exists an elliptic curve of conductor $pN$ where $p$ is prime, $N$ is relatively prime to $p$, and $X_0(N)$ has genus $0$, then the Steinberg representation of $\SL_2(\F_p)$ is Mordell-Weil in dimension $1$.
\end{thm}

\begin{proof}
Let $J_{p,N}$ denote the Jacobian variety of $X_{p,N}$.
We consider the diagram
\begin{equation*}
\begin{CD}
X_{p,N} @>>> J_{p,N} \\
@VVV  @VVV \\
X_0(N) @>>> J_{p,N}/\SL_2(\F_p). \\
\end{CD}
\end{equation*}
Let $A = J_{p,N}$ and $B = E$, where $E$ is any elliptic curve of conductor $Np$.
As $\Gamma_0(pN)\subset \Gamma_{p,N}$, $X_{p,N}$ maps onto $X_0(pN)$,
which maps onto $E$ by the modularity of elliptic curves over $\Q$.  Let $\pi\colon A\to B$
be a non-constant map factoring through the Jacobian variety of $X_0(pN)$.
Applying Proposition~\ref{realize}, it suffices to prove that $\Hom_{\C}(A,B)\otimes\Q$, regarded as a rational $\SL_2(\F_p)$ representation contains the Steinberg representation as a subrepresentation.

Let $P$ denote the image of $\Gamma_0(p)$ in $\SL_2(\F_p)$, i.e., the group of upper triangular matrices in $\SL_2(\F_p)$.  By construction,
$\pi$ is fixed by the action of $P$ on $A$.  Thus, the $\Q[\SL_2(\F_p)]$-submodule of
$\Hom_{\C}(A,B)\otimes\Q$ generated by $\pi$ is a quotient of $\Ind_P^{\SL_2(\F_p)}\Q$,
which is isomorphic to the direct sum of a trivial $1$-dimensional representation and the Steinberg representation.  We need only prove, therefore, that the action of 
of $\SL_2(\F_p)$ on $\pi$ is non-trivial, i.e., that $\pi$ does not factor through the maximal $\SL_2(\F_p)$-invariant quotient, $A_{\SL_2(\F_p)}$, of $A$.    However,
$$H^1_{\sing}(A_{\SL_2(\F_p)}(\C),\Q) \cong H^1_{\sing}(A(\C),\Q)_{\SL_2(\F_p)}
\cong H^1_{\sing}(X_{p,N}(\C),\Q)_{\SL_2(\F_p)},$$
which is trivial since the quotient of $X_{p,N}$ by $\SL_2(\F_p)$ is the genus $0$ curve
$X_0(N)$.

\end{proof}

\begin{cor}
The Steinberg representation of $\SL_2(\F_p)$ is Mordell-Weil for all primes $p<1000$.
\end{cor}

\begin{proof}
This follows immediately from the proposition by inspecting Cremona's tables \cite{Cr}.

\end{proof}

We remark that it is expected but not yet known 
that infinitely many primes satisfy the hypotheses 
of Theorem~\ref{steinberg}.  Work of Friedlander and Iwaniec \cite{FI}
gives some hope that this problem may be accessible.

We note that by restricting the Steinberg representation to the group $\Z/p\Z$ of unitriangular matrices in $\SL_2(\Z/\Z)$,
we deduce that the regular representation of $\Z/p\Z$ is Mordell-Weil in dimension $1$
for all primes $p<1000$.  In fact, one can prove more:

\begin{prop}
For every prime $p$, the regular representation of $\Z/p\Z$ is Mordell-Weil in dimension $1$.
\end{prop}

\begin{proof}
Mazur and Kur\v{c}anov observed \cite{Ku}, that under certain common conditions,
the rank of an elliptic curve over a $\Z_p$-extension of a number field is infinite.
For convenience, we use a more recent result due to Cornut \cite{Co} and Vatsal \cite{Va}.

We fix a non-CM elliptic curve $E$ with root number $-1$ and conductor $N<1000$.
We let $K_\infty$ denote the anti-cyclotomic $p$-extension of $\Q(i)$.
Every $p>1000$ is prime to $N$ and the conductor of $E$, and it follows 
that the rank of $E$ over $K_\infty$ is infinite.  Therefore there exists an abelian $p$-extension $K_{n+1}/K_n$ such that
$$\dim E(K_{n+1})\otimes\Q > \dim E(K_n)\otimes\Q > 0,$$
and it follows that the $\Gal(K_{n+1}/K_n)$-module $E(K_{n+1})\otimes\Q$
contains a copy of the regular representation of $\Z/p\Z$.

\end{proof}

\end{document}